\documentclass[12pt]{article}

\usepackage[utf8]{inputenc}

\usepackage{amssymb}
\usepackage{latexsym}
\usepackage{exscale}

\usepackage{tikz}
\usetikzlibrary{trees}

\usepackage{amssymb, amsmath}
\usepackage{enumerate}
\usepackage{enumerate}
\usepackage[top=1in,left=1in,right=1in,bottom=1in]{geometry}
         \usepackage{ulem}
          \usepackage{xspace}

  \newcommand{\jg}[1]{{\color{violet}{#1}}}
   \newcommand{\jb}[1]{{\color{blue}{#1}}}

\usepackage{verbatim}

\def \beq{\begin{equation}}
\def \eeq{\end{equation}}

\renewcommand{\rq}[1]{(\ref{#1})}

                                     \newtheorem{remark}{Remark}

\newtheorem{lemma}{Lemma}
\newtheorem{prop}{Proposition}
\newtheorem{thm}{Theorem}
\newtheorem{cor}{Corollary}

\newcommand{\bR}{{ \mathbb R  }}
\newcommand{\bC}{\Bbb C}
\newcommand{\bZ}{\Bbb Z}

\newcommand{\bN}{\mathbb{N}}
\newcommand{\bK}{\Bbb K}

\newcommand{\E}{\mathcal{E}}

\newcommand{\la}{\mbox{$\lambda$}}

\newcommand{\pa }{\partial }
\newcommand{\f}{\varphi}

\newcommand{\ep}{\epsilon}

\newcommand{\al}{\alpha }

\newcommand{\ga}{\gamma }

\newcommand{\La}{\Lambda }

\newcommand{\om}{{\omega}}

\def\<{\langle} \def\>{\rangle}

\newcommand\CM[1]{\vskip3mm\fbox{\parbox{5in}{#1}}\vskip3mm}
     
                \newcommand{\cN}{\mathcal{N}}
                
             \newcommand \tl {\widetilde}
             \newcommand\wrt{with respect to \xspace}
\title{ Null-controllability for the beam equation with structural damping.\\
Part 1. Distributed control}

\author{Sergei Avdonin, Julian Edward, and Sergei Ivanov }
\date{January 20, 2024}

\begin{document}
\maketitle

\begin{abstract}
Let $\Delta$ be the Dirichlet Laplacian on the interval $(0,\pi)$.
The null controllability properties of the equation
$$u_{tt}+\Delta^2 u+\rho (\Delta)^\al u_t=F(x,t)$$
are studied. Let $T>0$, and assume initial conditions $(u^0,u^1)\in Dom(\Delta)\times L^2(0,\pi)$. 
We first prove finite dimensional null control results: 
suppose
$F(x,t)=f^1(t)h^1(x)+f^2(t)h^2(x)$ with $h^1,h^2$ given functions. 
For $\al \in [0,3/2)$, we prove that there exist $h^1,h^2\in L^2(0,\pi)$ such that for any $(u^0,u^1)$, there exist $L^2$ null controls $(f^1,f^2).$ For
$\al< 1$ and $\rho <2$, we prove null controllability with $f^2=0$ and $h^1$  belonging to a large class of functions. For $\al\in  [3/2,2)$, we prove  spectral and null controllability both generally fail, but two dimensional weak controllability holds.
Our second set of results pertains to 
 $F(x,t)=\chi_\Omega(x)f(x,t)$, with $\Omega$ any open subset of $(0,\pi)$. For any $\al \in [0,3/2),$ we prove there exists a null control $f\in L^2(\Omega\times(0,T))$   To prove our main results, we use the Fourier method to rewrite the control problems as moment problems. These are then solved by constructing biorthogonal sets to the associated exponential families. These constructions seem to be non-standard and may be of independent interest.

\end{abstract}

\section{Introduction}%\label(Intro)

Let $\Delta$ be the  Laplacian: $\Delta =-\pa_x^2,$ with operator domain $H^2(0,\pi)\cap H^1_0(0,\pi)$. It is well known that this operator is self-adjoint with positive spectrum, and hence $\Delta^{\alpha}$ is defined by the Spectral Theorem for all $\alpha$.
We will study control problems for the equation
$$
u_{tt}+\Delta^2 u+\rho (\Delta)^\al u_t=Bf,\ x\in (0,\pi), \ t>0,
$$
with a positive constant $\rho$, $\al\in [0,2],$
and where $B$ is the control operator.

This system is actuated through a control mechanism prescribed by the operator $B$ (possibly unbounded to take into account trace operators prescribing
the boundary value of distributed states). Throughout this paper, controllability will always mean
the ability of steering any initial state $(u(x,0),u_t(x,0))$ to zero over a finite time by some appropriate
input function $f$ (i.e. exact controllability to zero or null controllability).

The term $(\Delta)^\al u_t$  models a specific dissipative effect, known as structural damping, when $\al \in (0,2)$. To the best of our knowledge, this was introduced in \cite{CR} assuming $\al =1$:
“The basic property of structural damping, which is said to be consistent with empirical studies,
is that the amplitudes of the normal modes of vibration are attenuated at rates which are proportional to the oscillation frequencies.” This model was also studied under the name “proportional
damping” (cf. \cite{Ba}). The quite different case $\al=2$ is known as “Kelvin–Voigt” damping. When
$B$ is the identity and $\al \in (0,2]$, this is the first class of parabolic-like control models considered
in \cite{LT},\cite{T2}, see also \cite{AL}.

The paper consists of two parts. In this first part we consider interior controllability, and in the second part  boundary controllability.

First, we consider two-dimensional interior control:
\begin{eqnarray}
u_{tt}+\Delta^2 u+\rho (\Delta)^\al u_t & = &f^1(t)h^1(x)+f^2(t)h^2(x) ,\nonumber \\
& &\ x\in (0,\pi), \ t>0,\label{beam3}\\
u(0,t)=u_{xx}(0,t)=u_{xx}(\pi,t) =u(\pi,t)& = & 0,\label{bc2inta}\\
u(x,0)=u^0(x),\ u_t(x,0)& = & u^1(x).\label{init2a}
\end{eqnarray}
Here $h^1,h^2$ will be fixed functions (profiles), and $f^1,f^2$ would serve as controls. 

Throughout this paper, we will denote $X^p=Dom(\Delta^{p/2})$, so $X^0=L^2(0,\pi)$, etc.
In the theorem below,  
\beq \label{Q}
Q(T)\leq \left \{
\begin{array}{cc}
C'/T, & \al \leq 1,\\
C'/T^{1/(3\al-2)}, & \al \in (1,3/2).
\end{array}\right .
\eeq
where the constant $C'$ depends on $\al$ and $\rho$ only. 
\begin{thm}\label{2d}
Consider the system \rq{beam3}, with boundary conditions \rq{bc2inta} and initial conditions \rq{init2a}. 
Suppose $(u^0,u^1)\in X^2\times X^0$ and $T>0$.

A) (One dimensional control.)
Suppose $\al = 0$, or  $\al \in (0,1]$ and $\rho \leq 2$. Set $h^2=0$. 
Suppose the Fourier coefficients 
$\{h^1_n\}$ of  $h^1\in L^2(0,\pi )$ satisfy $|h^1_n|\asymp 1/n^p$ for some  constant $p>1/2$.
{Then} there exists 
$f^1\in L^2(0,T)$ such that the solution $u$ to the system above solves
$$u(x,T)=u_t(x,T)=0,
$$
with 
$$
\| f\|_{L^2(0,T)} \leq Ce^{Q(T)}(\| u^0_{xx}\|_{L^2(0,\pi)}+\|u^1\|_{L^2(0,\pi)}).
$$
Here the constant $C$ is depends only on $\al,\rho$.

 B) (Two dimensional control.)
Suppose    $\al \in (0,1]$ and  $\rho > 2$, or $\al \in (1,3/2).$  Then there exist $h^1,h^2\in L^2(0,1)$ such that for any pair $(u^0,u^1)$,
  there exist 
$f^1,f^2\in L^2(0,T)$ such that the solution $u$ to the system above solves
$$u(x,T)=u_t(x,T)=0,
$$
with 
$$
\| f^1\|_{L^2(0,T)}+\| f^2\|_{L^2(0,T)}\leq Ce^{Q(T)}(\| u_{xx}^0\|_{L^2(0,\pi)}+\|u^1\|_{L^2(0,\pi)}).
$$ 
Here the constant $C$ is depends only on $\al$,$\rho$, and $h^1$,$h^2$.

C) If $\al \geq 3/2$, then the system is neither null controllable nor spectrally controllable with two dimensional (or with any finite dimensional) control.

D) (Two dimensional weak controllability.) Suppose $\al \in [3/2,2)$. Then there exist $h^1,h^2\in L^2(0,1)$ such that for any pair $(q^0,q^1)\in X^2\times X^0$ and any $\ep >0$,
  there exist 
$f^1,f^2\in L^2(0,T)$ such that the solution $u$ to the system above satisfies
$$\|u_{xx}(x,T)\|_{L^2(0,\pi )}+\|u_t(x,T)\|_{L^2(0,\pi)}<\ep .
$$
\end{thm}

\begin{remark}
Recall that spectral null-controllability is equivalent to the fact that 
jg{for} all the initial data $(u^0, u^1)$,
where $u^0$ and $u^1$ are eigenfunctions or zero, the system can be steered to rest and equilibrium. 
\end{remark}

To prove this result, we first apply the Fourier method, i.e. find the solution as  
$$
\sum a_k(t)e^{i\la_kt}\f_{|k|}(x),
$$
where $k$ runs over $\bK :=\bZ\backslash 0$, the sequence  $\La$ consists of the frequencies $\la_k$ arising in this method 
and $\f_k$ are the eigenfunctions of the Laplace operator. 
Then we rewrite the associated control problem as a moment problem. 

For part A, we need to distinguish the case $ \rho <2$ from all other cases.  In this case, the frequencies are separated. Hence
the associated exponential family, $\{ h_{|k|}^1e^{it\la_k},k\in \bK\}$, is minimal on $L^2(0,T)$, and hence the moment problem has a (formal) solution via the biorthogonal elements. The separation condition, together with the asymptotics of $\{ \la_k\}$, allow us to use a result in \cite{AIS} to  conclude the biorthogonal functions satisfy 
an exponential estimate that implies the formal solution converges in $L^2(0,T).$ 

For the case  $\rho=2$ and $\al=1$ the family has the form
$\E_1=\{e^{it\la_n},te^{it\la_n}: n\in \bN\}$.
We show that the family $\E_1$ is also minimal, and the elements of the biorthogonal family satisfy the necessary estimates to prove the theorem. 
The construction of the biorthogonal family, an adaptation of the argument of \cite{AIS}, is formulated in Proposition \ref{prop2} and might be of independent interest. In this context, we note that the construction of, and estimates on,  sets of functions
biorthogonal to $\{ e^{it\la_k}\}$
have been a subject of considerable research, see for instance \cite{BO} and references therein. 

For $\rho=2$ and $\al <1$, we have a single double frequency. This can be treated similarly to $\rho=2$ and $\al =1$.

For the part B, the frequency set $\{ \la_n\}$  no longer necessarily satisfies the separation condition. Moreover it is possible for some $\rho$ that two elements (no more!) can coincide and one dimensional control fails. Associated to the moment problem in this case is
the `vector' exponential family
\beq\label{cE} 
\E=\left \{ \binom{ h^1_{|k|}}{ h^2_{|k|}}e^{i\la_k t},k\in \bK\right \},
\eeq
where $h_k^j$ are the Fourier coefficients of the profiles $h^j$. 
This family $\E$ can be made  minimal by carefully choosing the Fourier coefficients of $h^1$, $h^2$.
Roughly speaking we split our vector exponential family into two 
 orthogonal families, and thus the original moment problem is split into two solvable moment problems, one for $f_1$ and one for $f_2$.

The proof of part C follows from the theory of moment problems, because the frequency set does not satisfy the Blaschke condition. The proof of part D uses an adaptation of the splitting argument from part B, together with the following sufficient condition for weak linear independence of a family of exponentials $\{ e^{\la_nt}; n\in \bN\}$ on $L^2(0,T)$, where we assume $\la_n >0$: 
$$
\lim_{n\to \infty}\frac{\ln (n)}{\la_n}=0.
$$
\

The next control problem is the following. 
Let $\Omega$ be an open subset of $(0,\pi).$
Suppose we have the initial boundary value problem (IBVP)
\begin{eqnarray}\label{beam}
u_{tt}+\Delta^2 u+\rho (\Delta)^\al u_t & = & \chi_{\Omega}f(t,x), \ x\in (0,\pi), \ t>0,\label{beam2}\\
u(0,t)=u_{xx}(0,t)=u_{xx}(\pi,t) =u(\pi,t)& = & 0\\
u(x,0)=u^0(x),\ u_t(x,0)& = & u^1(x).\label{init2}
\end{eqnarray}
Here $\chi_{\Omega}$ is the characteristic function of $\Omega$.

\begin{thm}\label{thm2}
Let $\al \in [0,3/2).$   Given 
$(u^0,u^1)\in X^2\times X^0$ and $T>0$, there exists $f\in L^2(\Omega\times (0,T))$ such that the solution $u$ to the system \rq{beam2}-\rq{init2} solves
$$u(x,T)=u_t(x,T)=0,
$$
with 
$$\| f\|_{L^2(\Omega\times (0,T))}\leq Ce^{Q(T)}(\| u^0_{xx}\|_{L^2(0,\pi)}+\|u^1\|_{L^2(0,\pi)}).
$$
Here $Q(T)$ is as in Theorem 1,  and the constant $C$ depends on $\al$ and $\rho$.
\end{thm}

For the  proof, we again reduce the control problem to the moment one, now  \wrt  the  exponential family
$$
\E_2=\jg{\{e^{i \la_k t}} \f_{|k|}(x)\big|_\Omega, k\in \bK\}.
$$
For the case of non-separated spectrum we split the the family using the fact that 
the angle between the eigenfunctions $\f_n$ and $\f_m$ in $L^2(a,b)$ is separated from zero.

%To prove this result, we \jg{apply the Fourier method and} 
%\sout{use the asymptotics of $\la_n^{\pm}$ to} prove the applicability of an Ingham-type inequality. For that % % we generalize Theorem 1 in \cite{AIS} to include the terms of the form $t e^{i\la_n t}$ and use the fact that % % the observability inequality is equivalent to the null-controllability. 

In the Part 2 of the paper, to be published separately, we will deal with  boundary controllability. We will study the dynamical system of the same form
and   discuss several approaches to treat non-homogeneous boundary conditions 
\begin{eqnarray}
u(0,t)=u_{xx}(0,t) & = & 0,\\
u(\pi,t) & = & f(t),\\
u_{xx}(\pi,t) & = & g(t).
\end{eqnarray}
We will then prove null controllability for this system.

This paper is organized as follows. In the next subsection, we compare our results with the literature. In Section 2.1, we discuss the spectral solution of the uncontrolled system, discussing how  $\al, \rho$ determine the separation properties of the frequencies {$\{ \la_n\}$}. In Section 2.2, we adapt an argument from \cite{AIS} to prove a existence of a biorthogonal to $\{e^{\la^{+}_n t},\,te^{\la^{+}_n t} \}$ family of functions, with the norms satisfying an exponential estimate, that will
be used to solve the moment problems associated to Theorems \ref{thm2} and \ref{2d}. 
 Theorem
\ref{2d} is proven in Section \ref{S4} and Theorem \ref{thm2}  is proven in Section  \ref{proofTh1} .

\subsection{Literature review}

We first compare our Theorem \ref{thm2} to the relevant literature. 
Lasiecka and Triggiani, in \cite{LT}, considered the case an abstract system which, for dimension 1, can be reduced to the beam equation with control distributed throughout the interval, and $\al \in [1,2)$.
Excepting the case $\al =1, \rho =2$, they prove null controllability. An important ingredient in their calculations is Parseval's Formula, which  requires $\Omega=(0,\pi).$
The well-posedness and regularity of the equation, and also for plate equations in higher dimensions,  is discussed in \cite{T1} using the theory of analytic semigroups, but assuming $\al$ is an integer.

Miller \cite{mil} considered  \rq{beam} in a bounded domain of $\bR^n,$ with distributed 
controls supported on a subset of the interior, which for $n=1$ could a arbitrary open subset $\Omega\subset (0,\pi )$. Miller proves null-controllability for $\al \in (1/2,3/2)$  when $\Omega$ is a proper subset, and for any $\al <1$ if $\Omega =(0,\pi ).$ In place of Parseval's Formula, he uses the inequality (which follows from a Carleman estimate due to Lebeau and Robiano)
$$\int_{\Omega}\left|\sum_{j\leq \om_j}c_j\phi_j(x)\right|^2dx\geq C_1e^{-C_2\om_j} \sum_{j\leq \om_j}|c_j|^2,
$$
along with the natural damping properties of the system. Here $\{ \om_j, \phi_j\}$ are the spectrum and corresponding orthonormal eigenfunctions for the Laplacian. Edward \cite{E} considered also interior control for $\al <1/2$, and for $\al =1/2$ with small $\rho$, and proved null controllability. 
More recently, Mitra \cite{mit} considered the case $\al =1$ on the interval with periodic boundary conditions. Using a Carleman estimate, null controllability is proven for controls supported on an open subset of $(0,1).$

Since our Theorem \ref{thm2} covers $\al \in [0,3/2)$ for all $\rho$, this theorem can be viewed as complementary to the results of Miller, Edward and Mitra. And compared with Edward's result, our Theorem \ref{thm2} has the advantage of the estimate on the cost blowup rate. 

As a final remark about Theorem \ref{thm2}, the case $\al \geq 3/2$, for $\Omega$ a proper subset of $(0,\pi)$, remains open. If $\Omega=(0,\pi)$, null controllabilty was proven in \cite{LT}.

%Miller also proves notes: ``the controllability cost of a system is not increased by taking its tensor product with a contraction semigroup". As a consequence, he can prove that, in the case $\al =1$, one has boundary control for a rectangle with control supported on a single edge. 

Regarding our Theorem \ref{2d}, we are unaware of any related papers for the structurally damped beam equation with finite dimensional interior control. Well-posedness and regularity for one dimensional control is discussed in \cite{T1}.

We conclude this section with a brief discussion 
on the function $Q(T)$,  which appears in our upper bound,
$Ce^{Q(T)}$, on the cost of control. 
The rate at which the cost of a control blows up as $T\to 0^+$ has been the subject of interest,  
motivated by problems in non-linear partial differential equations and stochastic differential equations. 
Recall we showed \jg{the estimate \rq{Q}}
%$$
%Q(T)\leq \left \{
%\begin{array}{cc}
%C'/T, & \al \leq 1,\\
%C'/T^{1/3\al-2}, & \al \in (1,3/2).
%\end{array}\right .
%$$
 for some constant $C'$,  which follows from the construction of $Q$ given in \cite{AIS}. 
Miller also estimates $Q(T)$, and comparing his Corollary 1
 with our Theorem \ref{thm2}, one sees that his result for interior control is sharper for the power of $T$.
The problem was also studied in \cite{AL}.
The results there are given in an abstract setting, but restricting to  the 1-dimensional beam equation, the control function is assumed to be distributed throughout the interior,
and under this strong hypothesis the control cost has the much smaller upper bound $C/T^{\beta}$, where $\beta$ is determined by $\rho$, $\al$.

\section{Frequency set and biorthogonal functions}\label{spectr}
\subsection{Frequency set}\label{2.1}
Consider the eigenvalue problem
$$
\Delta \f  = \la \f ,\
\f (0) =\f(\pi)= 0.\\
$$
Clearly an orthonormal  basis of eigenfunctions is $\{ \f_n; \ n\in\bN\} $, 
$\f_n(x)=\sqrt{\frac{2}{\pi}}\sin ( n x)$, with corresponding eigenvalues 
$n^2.$
Obviously the $(\Delta)^\al$ has the same eigenfunctions, with corresponding eigenvalues $n^{2\al}$.

Consider the IBVP , whose solution we will refer to as the ``free" wave,
\begin{eqnarray} 
w_{tt}+\Delta^2 w+\rho(\Delta)^\al w_t & = & 0,  \label{w_prime1} \\
w(0,t) =w(\pi ,t)=w_{xx}(0,t)=w_{xx}(\pi ,t)& = & 0, \\
w(x,0) & =& w^0(x),  \\
w_t(x,0) & = & w^1(x).  \label{w_prime4}
\end{eqnarray}

%Assuming $v,w$ are sufficient regular, integration by parts gives 
%$$
%0=\int_0^1[w_tv-wv_t]_0^T+\int_0^T\int_0^1 %w[f''(t)x+\rho f'(t)\Delta^\al x]dx\ dt+\rho %\int_0^1[w\Delta^\al v]_0^T\ dx.
%$$

Set $w=\sum_{n=1}^\infty a_n(t)\f_n(x)$. Then 
\rq{w_prime1} implies
$$
\sum (a_n''+\rho a_n'n^{2\al}+a_nn^4   )\f_n(x)=0,
$$
hence
\beq
a_n''+\rho n^{2\al}a_n'+n^4a_n=0,\ \ \forall n\in \bN.
\eeq
%Here $w^j_n=\sqrt{2}\int_0^1\sin (\om_n x)w_j(x)dx$.
Solving  $\la^2+\rho n^{2\al}\la+n^4=0$, we get
$$ 
\la=
\frac{-\rho n^{2\al}\pm \sqrt{\rho^2n^{4\al}-4n^4}}{2} =: \la_n^{\pm}.
$$
Thus, if $\la_n^+\neq \la_n^-$,
\beq
w(x,t)=\sum_1^\infty (c_n^+e^{\la_n^+ t}+c_n^-e^{\la_n^- t})\f_n(x),\label{du}
\eeq
with coefficients $c_n^\pm$ determined by the initial conditions:
\begin{eqnarray*}
c_n^+ + c_n^-=&w_n^0,\\
\la_n^+c_n^+ + \la_n^-c_n^-=&w_n^1,
\end{eqnarray*}
where $ w_n^0$ and  $w_n^1$ are the Fourier coefficients of $ w^0$ and $w^1$.
This gives the following expression
\begin{eqnarray}
c_n^+=&(w_n^1-\la_n^-w_n^0)/q_n,\label{cn+}\\
c_n^-=&(-w_n^1+\la_n^+w_n^0)/q_n,
\label{cn-}
\end{eqnarray}
where
\beq
q_n=\la_n^+-\la_n^-=\sqrt{\rho^2n^{4\al}-4n^4}.\label{qn}
\eeq

%\CM{I return to write $q_n$ here, because they will appear later. Ofcourse we  can write everywhere $\la_n^+-\la_n^-$}

If for some $n$ we have $\la_n^+=\la_n^-$, 
then we change the corresponding term in $w$ to
$$
(c_n^+e^{\la_n^+ t}+c_n^-te^{\la_n^+ t})\f_n(x),
$$
with 
\begin{eqnarray*}
c_n^+&=&w_n^0,\\
c_n^-&=&(w_n^1-\la_n^+w_n^0)/(\la_n^++1).
\end{eqnarray*}

We now examine the gap properties of 
$\{ \la_n^{\pm}\}$. 
 In what follows we will use the \textit{frequency set}
\beq \label{La}
\Lambda=\{\la_k\}_{k\in\bK},\ 
\bK=\bZ\backslash \{0\},\ 
		\la_k =\left \{
		\begin{array}{cc}
			-i\la_k^+, & k>0,\\
			-i\la_k^-, & k<0.
		\end{array}
		\right .
\eeq	
Thus it is easy to see $\Lambda\subset\mathbb C^+$.  We need (following \cite{AIS}) 
to introduce a function $\nu :[0,\infty)\mapsto [0,\infty)$
which describes the density of $\Lambda$ 
$$ 
\#\{ \la_n\in \Lambda\setminus \la_k  :  |\la_n -\la_k| < r\} \leq \nu (r), \forall k.
$$
We will require that $\Lambda$ satisfies 
\beq
\nu(r)=0, \ r<R_0,\label{separ}
\eeq
 for a positive R. This assumption is equivalent to 
 $$ \inf_{k\neq n}|\la_k-\la_n|>0,$$
which we well refer to as separability of $\Lambda.$

{\bf Example 1.} It is not hard to see that for $\la_n =sgn (n)|n|^p$, with $p > 1$ and $n\in \bN$, we have 
\beq \label{nu_asy} 
\nu(r)\asymp r^{1/p},
\eeq
for large $r$, {see also \cite{AIS}}.

Thus, if one assumes $\al \in [0,3/2)$, 
then the asymptotics below will  show that 
both sequences  $\{\la_k: k<0\}$ and  $\{\la_k:k>0\}$  satisfy \rq{separ}, \rq{nu_asy}
with $p=2$ for $\al\le 1$ and $p=2\al$ for $\al>1$.

\begin{lemma}\label{separat}

\ 

(1) For  $\al \in [0,1]$ and $\rho <2,$
 the frequency set is separable,

(2) For $(\al,\rho)=(1,2)$,   we have 
$$
\la_n^{+}= \la_n^{-}
$$ 
for all $n$. Also, for $\rho =2$ and any $\al\geq 0$, we have
$\la_1^+=\la_1^-.$

(3) For (i) $\al=1$ and $\rho>2$,  or (ii)  $\al>1$, there are an infinite number of $\rho$ such that the sequence $\Lambda$ contains two equal elements, i.e. for some $m$ and $n$
with $m\ne n$ we have
\beq \label{mn}
 \la_m^+= \la_n^-.
\eeq 

(4) For   $\al{\geq} 3/2$,  the frequency set does not satisfy the Blashke condition  \cite{Koosis}, in other words,
$$
\sum\left| \Im \frac1{i\la_n^+}\right| =\infty
$$

(5) For $\al =2$, we have $\la_n^+=O(1)$.

\end{lemma}
 \textbf{Proof:} 
(1)  (i) 
For $\al<1$ we have the following asymptotics: 
$$ 
\la_n^{\pm}= -\frac12\rho n^{2\al}(1+o(1))\pm in^2(1+o(1)).%\label{neq}
$$
From here we see that the frequency set is separable for large $n$.   Also for $\rho<2$
the set $\Lambda$ has no coinciding elements because the real parts of 
the branches $\la_n^{\pm}$ are strictly  increasing in $n$ while two branches have opposite signs of the imaginary parts.

(2) This is easily verified.

(3) (i) Let  $\al =1$, and $\rho > 2$.
Then 
$$
\la_n^+=n^2(-\rho /2+\sqrt{\rho^2/4 -1}), \la_n^-=n^2(-\rho /2-\sqrt{\rho^2/4-1}).
$$
Setting $r=(\rho+\sqrt{\rho^2-4})/2$, and we have  
\beq
\la_m^+= -\frac{2}{r} m^{2},\ \la_n^-= -\frac{r}{2}n^{2}.\label{bad}
\eeq
Take $m<n$. Then \rq{bad} implies that we can find $\rho$ satisfying  \rq{mn}. 

(ii) Let $\al \in (1,3/2)$. Then
\beq \label{bad1}
\la_n^+= -\frac{n^{4-2\al}}{\rho}\Big(1+o(1)\Big)\mbox{ and }\la_n^-= -\rho n^{2\al}\Big(1+o(1)\Big) .
\eeq

The situation is similar to 3(i). Fix $m$ and $n$. 
If $\rho$ runs $(0,\infty)$ the main terms of the branches in \rq{bad1} change from $-\infty$
to 0 and from 0 to $-\infty$, This means that we can find $\rho$ satisfying  \rq{mn}. 

(4) The statement can be checked directly.

(5) Follows immediately from \rq{bad1}.
 $\Box$

\begin{remark}
The case $\al =2$ is known as Kelvin-Voight damping, and because of part (5) of the lemma, the methods of this paper mostly cannot be used. In the case of controls distributed on $(0,\pi)$, this case is discussed in \cite{LT}, also see \cite{AI}.
\end{remark}

\begin{comment}
 hence the set $\{ e^{i\la_n^+ t}\}$ will not be minimal. 
\end{comment}

Finally, we discuss the properties of the solution to the IBVP \rq{w_prime1}--\rq{w_prime4}.
By \rq{du} and the asymptotics of $\{\la_n^{\pm}\}$, we have 
\begin{lemma}
The mapping $(w^0,w^1)\mapsto w$ is a continuous map 
$$X^2\times X^0\mapsto C(0,T;X^2\times X^0))\cap C^1(0,T; X^0).
$$ 
\end{lemma}
 {\textbf{Proof:} 
 Assume $(w^0,w^1)\in X^2\times X^0\times L^2(0,\pi )$, so that 
$$ 
\sum n^2(w_n^0)^2<\infty,\ \sum (w_n^1)^2<\infty.
$$
Now $w(\cdot,t)\in X^2\times X^0$ iff
$$
\sum n^2|c_n^+e^{\la_n^+ t}+c_n^-e^{\la_n^- t}|^2<\infty.
$$
Because all $\la_n^\pm$ have a negative real part, it is enough to check that 
$$
\sum n^2(|c_n^+|^2+|c_n^-|^2)<\infty.
$$
By \eqref{cn+},\rq{cn-}, and \rq{qn}, it suffices to show that 
\beq \label{conv}
\sum n^2(|(w_n^1)^2+(|\la_n^+|^2+|\la_n^-|^2)|w_n^0|^2)/q_n^2<\infty.
\eeq
It is easy to see that for $\al\in[0,2],\,$
$n/q_n$ and $\la_n^\pm/q_n$ are bounded and 
\eqref{conv} is correct.
Because the series converges uniformly in $t$, we have
$w\in C(0,T;X^2\times X^0 ).$

The rest of the lemma can be proved similarly.
$\Box$}

\subsection{Biorthogonal Functions}

An important part of our proof is to construct  suitable sets of biorthogonal functions associated to $\{ e^{\la_n^{\pm}t}:n\in \bN\}$.

For completeness, we begin by providing results from \cite{AIS}. Let ${\bf L}$ be a subset of $\bZ$.
\begin{thm}(\cite{AIS})\label{AIS}
Suppose $\Lambda =\{ \mu_l:l\in {\bf L} \} \subset \bC^+$. Let $R_0>0.$
 
A- Suppose there exists a function $\nu (r)$ such that for all $l\in {\bf L},$
$\Lambda$ satisfies 
$$
\#\{ \mu_l\in \Lambda\setminus \mu_m  :  |\mu_l -\mu_m| < r\} \leq \nu (r) 
$$
with $\nu(r) = 0$ for $r<R_0$, and $\nu (r)/r^2$ integrable.  Then for any $\delta > 0$, $T_0 > 0$, there exists a constant $\tilde{C} = \tilde{C}(\delta, T_0)$
such that for any sequence $\{ a_l\}$, we have
$$\Sigma_l  |a_le^{i\mu_l \delta} |^2 \leq  \tilde{C}
\int_{0}^{T_0}
| \Sigma_l a_le^{i\mu_l t}|^2dt.
$$

B- Suppose $T=\delta=T_0$. If $\nu (r)\asymp r^{1/p}$ for large $r$, then the constant $\tilde{C}$ satisfies
$$\tilde{C}\leq C_1\exp (C_2/T^{1/(p-1)}),
$$
with constants $C_1,C_2$ independent of $T$.
\end{thm}
In what follows, we will often have $\{ \mu_l;l\in {\bf L}\}=\{ \la_k;k\in\bK\}$ , with $\la_k$ given by \rq{La}. Assuming $\al <3/2$,
the following estimates  follow from Section 2.
There exist positive constants depending on $\rho,\al$ such that one can choose $\nu (r)$ satisfying
\beq \label{nu}
C_0r^{\varkappa}\leq \nu(r)\leq C_1r^{\varkappa},
\eeq 
where $\varkappa =1/2$ for $\al <1$ or $\al=1,\rho\leq 2$, 
$\varkappa =1/2\al$ for $\al \in (1,3/2).$

 Proposition \ref{prop2}, below, is a  generalization of
 a result proven in \cite{AIS}, which in turn generalizes a result in \cite{H}.  
Recall
$$\< f,g\>=\int_0^Tf(t)\overline{g(t)}dt,
$$
where the bar denotes complex conjugation.
%We wish to construct a set biorthogonal to 
%$$\{ e^{i\la_k t},te^{i\la_kt}: \ k\in {\bf N}\}. $$
\begin{prop}\label{prop2}
Let $T>0$. %\sout{Let $\Lambda$ satisfy \rq{nu},\rq{asymp}, \rq{not0}.}
Suppose there exists a function $\nu (r)$ satisfying  the estimates \rq{nu}  \rq{separ}.
 Then there exists a family of functions $\{  g_{m,j}(t); m\in {\bf L} ,\ j=1,2 \} $ in $L^2(0,T)$ satisfying 
$$ \< g_{m,1},te^{i\la_nt}\>=0,\  \< g_{m,1},e^{i\la_nt}\> =\delta_{m,n},\ \< g_{m,2},te^{i\la_nt}\>=\delta_{m,n},\  \< g_{m,2},e^{i\la_nt}\> =0 .
$$
Furthermore, there exist positive constants $C_2,C_3$ depending  only on $R_0,T,C_0,C_1$ such for $j=1,2,$
$$
\| g_{m,j}\|_{L^2(0,T)}\leq C_2 e^{C_3(\Im (\la_m))^{\varkappa}},%\label{bd2}
$$
where all constants are defined in \rq{nu} or Theorem \ref{AIS}.
\end{prop}
An immediate consequence of the proposition is: 
\begin{cor}
    \label{cor1}
Assume the hypotheses of Proposition \ref{prop2}.
Then there exists a family of functions $\{ g_j(t): j\in \bK \}$ in $L^2(0,T)$ satisfying 
$$ \< g_j,e^{i\la_nt}\> =\delta_{jn}.
$$
Furthermore, there exist positive constants $C_2,C_3$ depending  only on $R_0,T,C_0,C_1$ such
\beq
\| g_j\|_{L^2(0,T)}\leq C_2 e^{C_3(\Im (\la_j))^{\varkappa}}\label{bd1}.
\eeq 
\end{cor}
{\bf Proof of Proposition \ref{prop2}:} We adapt the construction used in \cite{AIS}. 
We define 
$$F_{j,1}(z)=\left (\prod_{k\in {\bf Z},k\neq j}\left(1-\left(\frac{z-\la_j}{\la_k-\la_j}\right)^2\right)\right)^2, j\in {\bN}.
$$
Then
$$F_{j,1}(\la_k )=\delta_{j,k}, F_{j,1}'(\la_k)=0,\ j,k\in \bN ,$$
and  by \rq{nu} and  \cite[Lemma 3]{AIS}, $F_{j,1}(z)$ is entire of exponential type zero with 
$$
        |F_{j,1}(\la_j+z)| \le e^{2\theta(|z|)}     \qquad\mbox{for }z\in\bC,
$$
where
$$
        \theta(s)=2\int_0^{\infty} \frac{\nu(r)}{r}\frac{s^2}{s^2+r^2}\,dr.
$$
Thus $\theta$ is a positive increasing function with 
\beq \theta (s)\asymp s^{\varkappa},\ s \to \infty.\label{theta2}
\eeq
We now define 
$$F_{j,2}(z)=(z-\la_j)F_{j,1}(z).
$$
Then for all $j,k\in {\bf N},$
we have 
$$F_{j,2}(\la_k)=0,F'_{j,2}(\la_k)=\delta_{j,k}.
$$
 Furthermore,
by increasing $\theta$ slightly, we have 
\beq \label{theta}
|F_{j,2}(\la_j+z)|\leq e^{2\theta(|z|)}, \ z\in \bC, \ j\in \bN .
\eeq
{Indeed, we can replace $\theta(s)$ by $\theta(s)+\log(s+1)-\log \min|\la_j|$}.
In what follows, we will employ this slightly larger $\theta$.
In what follows, it will be convenient to set $\al_j=\Re (\la_j)$, $\beta_j=\Im (\la_j)$, so $\beta_j\geq 0$.

By \cite[Theorem 2]{AIS} there exists 
 an entire function $P$ having  the following  properties

\noindent (i)  $ |P(z)|\leq 1 \mbox{ for } z\in \bC_+, \mbox{ and }P(0)=1$,

\noindent(ii)  $P(is)$ is real and positive for $s\geq 0$, and there exists a positive constant $C_4$ with
\beq \label{pl}
P(is)\geq e^{-C_4{s}^\varkappa},  \ s>0,
\eeq
\noindent(iii)
\beq \label{pb}
|P(s)|\leq e^{Q(T)}e^{-3\theta (|s|)}, \ s\in \bR,
\eeq
with $Q(T)$ a constant,

\noindent(iv) $P(z)e^{-izT/2}$ is of exponential type $T/2$.

Furthermore, we have
\begin{lemma}\label{P}
For $r\geq 0$, there exists $C_P>0$ such that
$$|P'(ir)|<C_P.$$
\end{lemma}
{\bf Proof:} We recall some facts from \cite{AIS}.
$$P(z)=\prod_{n=0}^{\infty}\frac{1}{2}(1+e^{2ia_nz}),
$$
with $\{ a_n\}$ a positive sequence satisfying $\sum_n a_n =\delta/2$ for some $\delta >0.$ 
Since $ P'(z)=P(z) (\log P(z))'$ and $|P(z)|\leq 1$ in the upper half plane,
$$|P'(ir)| \leq  |(\log P(ir))'|=\left| \left[\sum_n \log \left(\frac{1}{2}+\frac{1}{2}e^{-a_nr} \right) \right]' \right|  $$
$$ \leq \left| \sum_n \frac{-a_je^{-a_nr}}{1+e^{-a_nr}}  \right| \leq C_1  \sum_n a_n  \leq C_P.\ \Box$$

\ 

We now continue with the proof of the proposition.
Define, for $n=1,2$
\beq \label{Gj0}
G_{j,n}(z)=\frac{F_{j,n}(z)P(z-\al_j)}{P(i\beta_j)},\ j\in \bN .
\eeq
By \rq{theta} and \rq{pb}, $G_{j,n}(s)\in L^2(-\infty,\infty)$.
Furthermore,  for all $j,k\in \bN$, we have $G_{j,1}(\la_k)=\delta_{jk}$, and $G_{j,2}(\la_k)=0.$
By \rq{Gj0},
we have 
$$
G_{j,1}'(\la_k)=\frac{F'_{j,1}(\la_k)P(\la_k-\al_j)
+F_{j,1}(\la_k)P'(\la_k-\al_j)}{P(i\beta_j)}
=\delta_{jk}\frac{P'(i\beta_j)}{P(i\beta_j)}
$$
and
$$
G_{j,2}'(\la_k)=\frac{F'_{j,2}(\la_k)P(\la_k-\al_j)
+F_{j,2}(\la_k)P'(\la_k-\al_j)}{P(i\beta_j)}
=\delta_{jk}.
$$

Define 
$$g_{j,2}(t)=
\frac{1}{2\pi}\int_{\bR}\overline{G_{j,2}(s)}e^{ist}ds
$$
and 
$$g_{j,1}(t)=
\frac{1}{2\pi}\int_{\bR}\overline{(G_{j,1}(s)-\frac{P'(i\beta_j)}{P(i\beta_j)}G_{j,2}(s))}e^{ist}ds
$$
with $j\in \bN$.

Then,  
$$ 
\< g_{j,1}(t),e^{i\la_k t}\>_{L^2(0,T )}=\delta_{j,k}, \ \< g_{j,2}(t),e^{i\la_k t}\>_{L^2(0,T)}=0,\  \forall k\in \bN.
$$
Also,
$$
\< g_{j,1}(t),te^{i\la_kt}\> =-i\frac{d}{d\la}\< g_{j,1}(t),e^{i\la t}\>|_{\la=\la_k}= 
-i(\overline{G_{j,1}'(\la_k)-\frac{P'(i\beta_j)}{P(i\beta_j)}G'_{j,2}(\la_k))}=0, \  n=1,
$$
$$
\< g_{j,2}(t),te^{i\la_kt}\> =-i\frac{d}{d\la}\< g_{j,2}(t),e^{i\la t}\>|_{\la=\la_k}= 
-i\frac{d}{d\la}\overline{G_{j,n}(\la)}|_{\la =\la_k}=\delta_{jk},  n=2.
$$

Furthermore, $F_{j,n}(z)$ is entire of  exponential type zero, and $e^{-izT/2}P(z-\beta_j)$ is entire of {exponential type  $T/2$ in both halfspaces},
and so by the Paley-Wiener Theorem, $g_{j,n}\in L^2(0,T)$.

We  see
$\{  g_{j,n}, \ j\in \bN , \ n=1,2\}$ is a  biorthogonal set to $\{  e^{i\la_kt},te^{i\la_kt}; k\in \bN \}.$
We now estimate the elements of this set. 
By \rq{theta},\rq{theta2}, \rq{pl}, and \rq{pb}, we have for $s\in \bR$
\begin{eqnarray*}
|G_{j,2}(s+\al_j)| & = & |F_{j,2}(\la_j+s-i\beta_j)P(s)/P(i\beta_j)|\\
&\leq &e^{2\theta(|s-i\beta_j|)}e^{Q(T)-3\theta(|s|)}/e^{-C_4(\beta_j)^\varkappa}\\
& \leq & e^{Q(T)-C_5{|s|^\varkappa}+(1+C_4)(\beta_j)^\varkappa}.
\end{eqnarray*} 
Since  the Fourier transform is unitary,
$$
%\| g_j\| \leq Ce^{(\beta_j)^\varkappa},\ 
\| g_{j,2}\| \leq Ce^{(1+C_4)(\beta_j)^\varkappa} , 
$$
with with $C$  depending only on $T,R_0,\ep, C_j,j=1-5.$  
Similarly, we estimate $g_{j,1}$ where we must use Lemma \ref{P}:
$$
\| g_{j,1}\| =
\| (G_{j,1}(s)-\frac{P'(i\beta_j)}{P(i\beta_j)}G_{j,2}(s))\|\leq  Ce^{(1+C_4)\beta_j^\varkappa}+ C_PCe^{(1+2C_4)\beta_j^\varkappa}.
$$
$\Box$

%In the notation of AIS-old-version, because
%$\la_n=\al_n+i\beta_n$ with $\al_n\asymp n^2$ and 
%$$\tilde{h}(r)\sim |r|^{1/2}, $$
%with $C_5$ independent of $j$.

\section{Proof of Theorem \ref{2d}}\label{S4}

Fix $h^1,h^2\in L^2(0,\pi)$. We consider
the following initial boundary value problem on $(0,\pi)\times (0,T)$
\begin{eqnarray}
u_{tt}+\Delta^2 u+\rho (\Delta)^\al u_t & = & h^1(x)f^1(t)+h^2(x)f^2(t),\label{z1}\\
u(0,t)=u(\pi ,t)=u_{xx}(0,t)=u_{xx}(\pi ,t) & = & 0,\\
u(x,0)=u^0(x),\ u_t(x,0)& = & u^1(x).\label{z4}
\end{eqnarray}
Here $(u^0,u_1)\in X^2\times X^0$. 
We wish to prove null-controllability.

We can represent the solution to \rq{z1}-\rq{z4} as a sum of a ``free" wave, corresponding to $f^1=f^2=0$, and a ``controlled" wave,
corresponding to $u^0=u^1=0$.
%We set
%$
%\ =
%\al_n+\beta_n, \ \la_n^-=
%-\al_n+\beta_n.
%$
Let us express the free wave, $u^{free}$, as a Fourier series. 
Suppose for $j=0,1$, the initial conditions have Fourier coefficients $\{ u^0_{n}\},\{ u^1_{n}\}$ respectively. Then, similarly to \rq{du}-\rq{qn}, if we assume $\la_n^+\neq \la_n^-$ for all $n$,
$$u^{free}(x,t)=\sum (c_n^+e^{\la_n^+t}+c_n^-e^{\la_n^-t})\f_n(x),$$
with
$$c_n^+=\frac{\la_n^-u^0_{n}-u^1_{n}}{\la_n^--\la_n^+},\ c_n^-=\frac{\la_n^+u^0_{n}-u^1_{n}}{\la_n^+-\la_n^-}.
$$
In the calculations below, we will assume $\la_n^+\neq \la_n^-$.
In the cases where $\la_n^+=\la_n^-$, the calculations below can be adapted by replacing $e^{\la_n^+t},e^{\la_n^-t}$ by $e^{\la_n^+t},te^{\la_n^+t}$. The details of the adaptation are left to the reader, but also see the paragraph at the end of the proof of part A below. 
Thus 
\beq
u^{free}(x,T)=\sum \ga_n^1\f_n(x):=\sum (c_n^+e^{\la_n^+T}+c_n^-e^{\la_n^-T})\f_n(x),\label{tf1}
\eeq
\beq\label{tf2}
u_t^{free}(x,T)=\sum \ga_n^2\f_n(x):=\sum (\la_n^+c_n^+e^{\la_n^+T}+\la_n^-c_n^-e^{\la_n^-T})\f_n(x).
\eeq
We now derive a formula for the controlled wave, denoted $u^f$ with $f=(f^1,f^2)$, and setting $u^0=u^1=0$. Let $h^j_n$ are the Fourier coefficients of $h^j$, $h^j_n=\< h^j,\phi_n\>$,
and
let $u^f(x,t)=\sum a_n(t)\f_n(x)$.
Putting this into \rq{z1}, we get the following family of ODE:
$$
a_n''+\rho n^{2\al}a_n'+n^4a_n= f^1(s)h^1_n+ f^2(s)h^2_n,\ a_n(0)=a_n'(t)=0,\ \forall n\in \bN.
$$
%Let
%\beq
%q_n=\la_n^+-\la_n^-=\sqrt{\rho^2n^{4\al}-4n^4}.\label{qn}
% \eeq
Then the solution to the ODE above is 
$$
a_n(t)=-\frac{1}{q_n}
\int_0^t\big ( f^1(s)h^1_n+ f^2(s)h^2_n\big )
\left(e^{\la^+_n(t-s)}-e^{\la^-_n(t-s)}\right)ds,\ n\in \bN.
$$
Comparing this with \rq{tf2}, we see that null controllability in time $T$ is equivalent to 
$$
\ga_n^1=\frac{1}{q_n}
\int_0^T\big ( f^1(s)h^1_n+ f^2(s)h^2_n\big )
\left(e^{\la^+_n(T-s)}-e^{\la^-_n(T-s)}\right)ds,
$$
$$
\ga_n^2=\frac{1}{q_n}
\int_0^T\big ( f^1(s)h^1_n+ f^2(s)h^2_n\big )
\left(\la_n^+e^{\la^+_n(T-s)}-\la_n^-e^{\la^-_n(T-s)}\right)ds,
$$
or, equivalently,
\beq \label{mo3}
\zeta_n^1=:{q_n}\frac{\ga_n^1-\frac{\ga_n^2}{\la_n^+}}{(-1+\frac{\la_n^-}{\la_n^+})}=
\int_0^T\big ( f^1(s)h^1_n+ f^2(s)h^2_n\big )
e^{\la^-_n(T-s)}ds,\ n\in \bN,
\eeq
\beq \label{mo4}
\zeta_n^2=:{q_n}
\frac{\ga_n^1-\frac{\ga_n^2}{\la_n^-}}{(1-\frac{\la_n^+}{\la_n^-})}=
\int_0^T\big ( f^1(s)h^1_n+ f^2(s)h^2_n\big )
e^{\la^+_n(T-s)}ds,\ n\in \bN.
\eeq
We set $\zeta_k=\zeta_k^1$ for $k>0$, $\zeta_k=\zeta_{-k}^2$ for $k<0$.
For $j=1,2$, we extend $h^j_n$ to $\bK$ by $h^j_{k}:=h^j_{|k|}$.
Recall
$$
\la_k =\left \{
		\begin{array}{cc}
			-i\la_k^+, & k>0,\\
			-i\la_k^-, & k<0.
		\end{array}
  \right .
$$
Hence
 system \rq{mo3},\rq{mo4} can be rewritten in terms of $i\la_n$:
\beq \label{mo51}
\zeta_k=
\int_0^T\big ( f^1(s)h^1_k+ f^2(s)h^2_k\big )
e^{\jg i\la_k(T-s)}ds,\ k\in \bK.
\eeq

\begin{remark}
Relating the above formula for our discussion in the introduction leading to \rq{cE},
 we have obtained the moment problem in $L^2\big((0,T)\times \bC^2\big)$
\wrt the solution $\tilde F(t)=\binom{f^1(T-t)}{f^2(T-t)}$ and the exponential family complex conjugate to \rq{cE}.
\end{remark}

{\bf {Proof of part A}}.

In this case, set $f^2=0.$ Assume for the moment $0<\al<1,\rho<2$. The other cases will be addressed at the end of the paragraph.
We have
$$\la_n^-= -in^2(1+o(1))-\frac{\rho n^{2\al}}{2}(1+o(1)),\ 
\la_n^+= in^2(1+o(1))-\frac{\rho n^{2\al}}{2}(1+o(1)).
$$
We see that the set $\Lambda$
satisfies the hypotheses of Corollary \ref{cor1}, and hence $\{ e^{i\la_kt},k\in \bK\}$ admits a biorthogonal family $\{ g_k,k\in \bK\}$ satisfying \rq{bd1}. It follows from \rq{mo51} that 
$$f^1(t):=\sum_{j\in \bK}\frac{\zeta_j}{h^1_j}\bar g_j(T-t) $$ 
formally satisfies  the moment problem. 
Recall the hypothesis that $|h^1_k|\asymp |k|^{-p}$ for some  positive constant $p$.
Combining  this with \rq{bd1} (with $\kappa =1/2$), \rq{tf1}, \rq{tf2}, \rq{mo3}, \rq{mo4},  we get 
$f^1\in L^2(0,T)$.

We now discuss the case  $\al=1$. 
If $\rho<2$
the frequency set satisfies the hypotheses of Corollary \ref{cor1}, and we can argue the same as the case $\al<1,\rho<2.$ 
If $\rho=2$, then we have $\la_n^+=\la_n^-$ for all $n\in \bN$. In this case, we can represent the control problem as a moment problem using the family $\{ e^{t\la_n^+},te^{t\la_n^+}: n\in\bN\}$. Then by Proposition \ref{prop2}, there exists a biorthogonal family of functions which can be used to solve the moment problem. The details are left to the reader.

Finally, suppose $\al=0$. A simple calculation shows multiple frequencies only arise
when $\rho=2n$. For $\rho \neq 2n$, we can argue as in the case $\al\in (0,1),\rho <2$ to prove the theorem. For $\rho=2n_0$ for some $n_0\in \bN$,  
we have $\la^+_n=\la^-_m$ if and only if $m=n=n_0$, and in this case we can still apply Proposition \ref{prop2} to obtain a biorthogonal family of functions which can be used to solve the moment problem.

\begin{remark} If $\rho>2$ and $\al \leq 1$,   the frequency set will have multiplicities for various values of $\rho$ and $\al$, in which case  Proposition \ref{prop2} won't apply, so we are in the situation Part B. Indeed, the sets $\{\la_k\}_{k<0}$ and  $\{\la_k\}_{k>0}$ separately 
satisfy the Proposition  \ref{prop2}, and the clusters (if any!) consist of two points.
\end{remark}

{\bf {Proof of part B}} 

We will present the proof for $\al \in (1,3/2)$; it will be easy to see that the case $\al \leq 1, \rho >2$ can be covered by the same argument. 

Recall  that, for $\al>1$, 
$$
\la_n^-\asymp -\rho n^{2\al},\ \la_n^+\asymp -\frac{1}{\rho}n^{4-2\al}.
$$
The difficulty in solving the moment problem, \rq{mo51}, is that we do not know whether 
the sequence
$\{ \la_k: k>0\}$ is  separated from
$\{ \la_k: k<0\}.$ We address this as follows.

For $\ep>0$, we will refer to the pair $(\la_n^+,\la_l^-)$ as an $\ep$ cluster if
$|\la_n^+- \la_l^-|<\ep$.
Let $\ep >0$ be sufficiently small that any $\ep$ clusters involving element of $\{ \la_n^+\},$
$\{ \la_l^-\}$ will involve only two elements. Let $\iota$ be the bijection within the set of $\ep$ clusters that maps  $\la^+_n$ to its cluster-counterpart $\la^-_l$.
Let 
\beq \label{cN}
{\cal N}^+=\{ n\in \bN: \exists l=\iota (n) \mbox{ such that } 
|\la_n^+- \la_l^-|<\ep \}=\mbox{domain}(\iota),\mbox{ and } {\cal N}^-=\mbox{range}(\iota).
\eeq
It is worth noting that if ${\cal N}^+\cap {\cal N}^-=\emptyset$, then the construction of $\{ h^1_n,h^2_n\}$ is easy: it suffices to define
\beq \label{Bmo1}
h^1_n=1/n, h^2_n=0\mbox{ if }n\in \bN\setminus {\cal N}^-,\mbox{ and } h^1_n=0, h^2_n=1/n\mbox{ if }n\in {\cal N}^-.
\eeq
Indeed,  the moment equalities \rq{mo3} and \rq{mo4} then take the form of the  moment equalities \wrt two separated sets of  exponentials 
\beq \label{Bm2}
\xi_n^1=(\frac 1n e^{\la^-_n t},\tl f_1),\ n \in  \bN\setminus {\cN}^+
\eeq 
\beq \label{Bm3}
\xi_n^2=(\frac 1n e^{\la^+_n t},\tl f_2),\ n \in {\cN}^+.
\eeq 
Here $\tl f_j(t) = \bar f_j(T-t)$, and the set $\{\xi_n^j\}$ is the renumbered set $\{\zeta_n^j\}$

Let us consider the general case: $\cN^\pm:=\cN^+\cap \cN^-$ is not empty.  The goal is to obtain also in this case two moment problems \wrt   
separated sets of exponentials. The first step is the first step in \rq{Bmo1}:
$$
h^1_n=1/n, h^2_n=0\mbox{ if }n\in \bN\setminus {\cN}^+.
$$
At this moment we have the moment equalities for \rq{Bm2} or for a part of \rq{mo3}. Evidently,
the set $\{\la^-_n \}_{n\in \bN\setminus {\cN}^+}$ is separated.

The second step is close to the second step in  \rq{Bmo1}:
$$
h^1_n=0, h^2_n=1/n\mbox{ if }n\in {\cN}^+\setminus \cN^\pm.
$$
The corresponding moment equalities are the part $\cN^+\setminus {\cN}^\pm$
of \rq{Bm3}.  The set $\{\la^+_n \}_{n\in \cN^+\setminus {\cN}^\pm}$ is separated.

Because of the asymptotics of $\la_n^{\pm}$, there exists $M>0$ such that $n>M$ implies
\beq \label{io1}
n>\iota(n).
\eeq
In the calculations that follow, we will assume \rq{io1} holds for  all $n$, leaving the simple adaptations for the general case to the reader.

The third step will involve an induction in which we define $h^j_n$ for $n\in \cN^{\pm}$.
 Let $m$  the smallest element in $\cN^\pm$, and let $l=\iota(m).$ 
 Because  $l<m$, we have 
 $$l\in \cN^-\setminus \cN^{\pm}\subset \bN \setminus \cN^+.$$
Hence $h^1_l=0,h^2_l=1/l$.
Thus we choose 
$$
h^1_m=1/m,h^2_m=0.
$$
For $n\in \cN^{\pm}$ for $n>m$, we carry out the following inductive step: 
\beq \label{Bmo2}
\mbox{ if } 
h^2_{\iota (n)}=0\mbox{  then set  } h^1_n=0, h^2_n=1/n,
\eeq
\beq \label{Bmo3}
\mbox{and if }
h^1_{\iota (n)}=0\mbox{  then set  } h^1_n=1/n, h^2_n=0.
\eeq
In the first case, the corresponding moment equality has the same form as in \rq{Bm3}, and in the second case,
 the moment equality has the same form as in \rq{Bm2}. The key point is that we never have $\la^+_n,\la^-_{\iota(n)}$ both appearing in one of \rq{Bm2},\rq{Bm3}. 
 In the other words, 
similarly to \rq{Bm2}, \rq{Bm3}
we obtain  moment problem \wrt  $f_1$ and $f_2$ with  two sets of "scalar " exponentials. Moreover, each set is a separated subset of
all $(\la_n^+,\la_l^-)$ and satisfy the Proposition 1.

\

{\bf Proof of part C }

Here, we will use the theory of the moment problem and its application to control problems, (\cite[Ch. I, V]{AI}).
If $\al\geq 3/2$, then by Lemma \ref{separat} the sequence $\{ \la_n^+\}$ fails the Blaschke condition,
and hence the exponential family is not minimal.
Thus both null controllability and spectral controllability will fail. 

\ 

{\bf Proof of part D}

We now prove the weak controllability for $\al\in [3/2,2).$
Recall 
$$
\la_n^+= -\frac{n^{4-2\al}}{\rho}\Big(1+o(1)\Big)\mbox{ and }\la_n^-= -\rho n^{2\al}\Big(1+o(1)\Big) .
$$
It is easy to see that the sets $\{{\la^+_n}\}$ and $\{{\la^-_n}\}$
 are each simple, and 
$$
\lim_{n\to \infty}\frac{\ln (n)}{\la_n^{\pm}}=0.
$$
Hence, by  (\cite{AI}, Theorem II.6.3 and ), the families
$\{ e^{\la_n^{+}t}\},\{ e^{\la_n^{-}t}\} $ are each weakly linearly independent.
But their union might not be, due to multiple frequencies. For this reason, we need to use two dimensional control. 

 We will present the proof for $\rho>2,$ in which case 
 $ \la^+_n> \la^-_n$ for all $n$. The adaptations for  cases $\rho=2$ and $\rho <2$ will be indicated at the end of the proof.
 
We express the controllability problem using moment problems \rq{mo3},\rq{mo4}:
$$
\zeta_n^1=
\int_0^T\big ( f^1(s)h^1_n+ f^2(s)h^2_n\big )
e^{\la^-_n(T-s)}ds,\ n\in \bN,
$$
$$
\zeta_n^2=
\int_0^T\big ( f^1(s)h^1_n+ f^2(s)h^2_n\big )
e^{\la^+_n(T-s)}ds,\ n\in \bN.
$$
We now adapt the construction of $h^1,h^2$ from part B. 
Recall that any multiple frequency will have multiplicity at most 2.  Let $\iota$ be the bijection within the set of double frequencies that maps  $\la^+_n$ to its counterpart $\la^-_l$.
Because of the asymptotics of $\la_n^{\pm}$, there exists $M>0$ such that $n>M$ implies
\beq \label{io2}
n>\iota (n).
\eeq 
In the calculations that follow, we will assume \rq{io2} holds for  all $n$, leaving the simple adaptations for the general case to the reader.
Let 
$${\cal N}^+=\{ n\in \bN: \exists l=\iota (n) \mbox{ such that } 
\la_n^+= \la_l^-\}=\mbox{domain}(\iota),\mbox{ and } {\cal N}^-=\mbox{range}(\iota).
$$
We can now argue exactly as in the proof for two dimensional null-controllability for $\al\in (1,3/2)$
to construct, from ${\cal N}^{\pm}$, the  functions $h^1,h^2.$ As in part B, the moment problem turns into two distinct moment problems, one each for $f^1$ and $f^2$, and in each the frequency set is simple, so the associated exponential families are weakly linearly independent. By (\cite{AI}, Thm. III.3.3), weak controllability follows.

 Finally, if $\rho \leq 2$, than it is possible that there exists $m$ such that $\la_m^+=\la_m^-$, and also possibly a finite number of complex frequencies.  Thus we must extend our exponential family to
$\{e^{\la_n^+t},te^{\la_m^+t},n\in \bN\},$ with a finite number of distinct non-real frequencies.
A careful reading of the proof of
(\cite{AI}, Theorem II.6.3) shows that this extended family remains weakly linearly independent. The remaining adaptations of the proof of above are left to the reader.

Our proof of weak controllability is complete.$\Box$

\section{Interior control of structurally damped beam. 
The proof of Theorem \ref{thm2}} \label{proofTh1}

Let $\Omega \subset (0,\pi )$ be a proper open subset. We consider:
\begin{eqnarray}
u_{tt}+\Delta^2 u+\rho (\Delta)^\al u_t & = & \chi_{\Omega}(x)f(x,t), \ x\in (0,\pi ), \ t>0,\label{z1i}\\
u(0,t)=u(\pi ,t)=u_{xx}(0,t)=u_{xx}(\pi ,t) & = & 0,\\
u(x,0)=u^0(x),\ u_t(x,0)& = & u^1(x).\label{z4i}
\end{eqnarray}
Here $u^0\in H^2\cap H_0^1,\ u^1\in L^2$.

%\sout{Since our goal is to prove null controllability, we can assume, 
%Without loss of generality,} { Evidently that it is enough to prove the theorem 
%for the case   $\Omega =(a,b)\subset (0,\pi)$.

First, we represent the solution to \rq{z1i}-\rq{z4i} 
using the calculations and notation from the previous section, see \rq{tf1},\rq{tf2}. 
In this section, we will assume $\la_n^+\neq \la_n^-$ for all $n$, leaving the simple adaptations in the other case to the reader.
Thus at $t=T$, we have 
the free wave satisfying 
$$
u^{free}(x,T)=\sum \ga_n^1\f_n(x):=\sum (c_n^+e^{\la_n^+T}+c_n^-e^{\la_n^-T})\f_n(x),
$$
\beq\label{tf2i}
u_t^{free}(x,T)=\sum \ga_n^2\f_n(x):=\sum (\la_n^+c_n^+e^{\la_n^+T}+\la_n^-c_n^-e^{\la_n^-T})\f_n(x).
\eeq
We now adapt the argument of the previous section to derive
 a formula for the controlled wave, denoted $u^f$. Let $f_n(t)$ are the Fourier coefficients of $\chi_{\Omega}(x)f(x,t)$, so
$$
f_n(t)=\int_0^{\pi}\chi_\Omega(x)f(x,t)\phi_n(x)\ dx=\int_\Omega f(x,t)\phi_n(x)\ dx,
$$
and
let $u^f(x,t)=\sum a_n(t)\f_n(x)$.
We have the set of ODEs
$$
a_n''+\rho n^{2\al}a_n'+n^4a_n=f_n(t),\ a_n(0)=a_n'(t)=0,\ \forall n\in \bN.
$$
Hence
$$
a_n(t)=-\frac{1}{q_n}
\int_0^tf_n(s)
\left(e^{\la^+_n(t-s)}-e^{\la^-_n(t-s)}\right)ds,\ n\in \bN.
$$
Comparing this with \rq{tf2i}, we see that null controllability in time $T$ is equivalent to 
$$
\ga_n^1=\frac{1}{q_n}
\int_0^Tf_n(s)
\left(e^{\la^+_n(T-s)}-e^{\la^-_n(T-s)}\right)ds,
$$
$$
\ga_n^2=\frac{1}{q_n}
\int_0^Tf_n(s)
\left(\la_n^+e^{\la^+_n(T-s)}-\la_n^-e^{\la^-_n(T-s)}\right)ds,
$$
or, equivalently,
\beq \label{mo3i}
\zeta_n^1=:{q_n}\frac{\ga_n^1-\frac{\ga_n^2}{\la_n^+}}{(-1+\frac{\la_n^-}{\la_n^+})}=
\int_0^T f_n(s)
e^{\la^-_n(T-s)}ds,\ n\in \bN,
\eeq
\beq \label{mo4i}
\zeta_n^2=:{q_n}
\frac{\ga_n^1-\frac{\ga_n^2}{\la_n^-}}{(1-\frac{\la_n^+}{\la_n^-})}=
\int_0^T f_n(s)
e^{\la^+_n(T-s)}ds,\ n\in \bN.
\eeq
We set $\zeta_k=\zeta_k^1$ for $k>0$, $\zeta_k=\zeta_{-k}^2$ for $k<0.$ Recall we have $\la_k=-i\la_k^+$ for $k>0$, and 
$\la_k=-i\la_{|k|}^-$ for $k<0.$
We extend $f_n$ to $\bK$ by $f_{k}:=f_{|k|}$, and similarly 
$\f_k(x)=\f_{|k|}.$
Then system \rq{mo3i}, \rq{mo4i} can be rewritten
\beq \label{mo51i}
\zeta_k=
\int_0^T\int_{\Omega}f(x,t)\f_k(x)
e^{\jg i\la_k(T-s)}\ dxds,\ k\in \bK.
\eeq

The remainder of this section will be devoted to solving this moment problem on $L^2(\Omega\times (0,T))$ by constructing a suitably bounded biorthogonal set to $\{\f_k(x)
e^{\jg i\la_k(T-s)},k\in \bK\} $. 
We assume first that for all $n$, $\la_n^+\neq \la_n^-$, which is equivalent to $\rho /2\neq n^{2-2\al}$. At the section's end, we briefly discuss the adaptations necessary in the other case. 

%Denote by $A(n,m)$  .  
\begin{lemma}\label{trig} The infimum of the angles between  $\f_n$ and $\f_m$ in $L^2(a,b)$ is positive.
% $$\inf_{m\ne n} A(n,m)>0.$$
\end{lemma}
The elementary proof of this lemma is deferred to the appendix. 

Let $\phi_n$ be the restrictions of the eigenfunctions $\f_n$ to $(a,b)$, normalized in the space $L^2(a,b).$
\begin{prop}\label{gen}
Let $T>0.$
Suppose $\al \in [0,3/2)$. 
Assume $\la_n^+\neq \la_n^-$ for all $n.$
Then there exists a set 
$\{ h_k(x,t):k\in \bK\}$ biorthogonal to 
$\{ \phi_{k}(x)\exp (i\la_kt),k\in \bK \}$ in $L^2((a,b)\times (0,T))$. Furthermore, there exist positive constants $C_2,C_3$, depending, on  $a,b,T$, such that 
\beq \label{biest1}
\int_0^T\int_a^b|h_k(x,t)|^2dxdt\leq C_2 \exp(C_3(Im{\la_k})^{\varkappa}).
\eeq
Here $\varkappa =1/2$ for $\al <1$ or $\al=1,\rho\leq 2$, and $\varkappa =3-2\al$ for $\al \in (1,3/2).$

\end{prop}
{\bf Proof:} 
First, if we assume $\al \leq 1$ and $\rho <2$, then
the frequency set is separated, so
by Proposition \ref{prop2} we can use $h_k(x,t)=g_k(t)\phi_k(x).$
Next, note $\rho=2$ is ruled out because we assume $\la_1^+\neq \la_1^-$.
In the remainder of the proof, we consider the harder cases $\al \in (1,3/2)$, or $\al \leq 1$ and $\rho>2$,
so that the union $\{ \la_n^+\} \cup \{ \la_n^-\}$ is not necessarily separated.

Recall that the cardinality of any cluster of frequencies can be at most two.  
We  use the notation  introduced in the previous section, in the proof of Theorem \ref{2d}. In particular, there exists $\ep >0$ such that 
 the set $\cN_\ep$  of $\ep$-close frequencies can be parametrized as 
$$
\cN_\ep= \{(l_n,m_n)\}_{n\in \bN}.
$$
Thus $(l_n,m_n)\in \cN_\ep$ means $|\la_{l_n}-\la_{m_n}|<\ep$.
We introduce now two reduced sets of the indices: 
$\bK_{\ep1}:=\bK \setminus \{m_n\},$ where  $\{m_n\}$ is the set of the second indices of  pairs from $\cN_\ep,$ and  $\bK_{\ep}:=\bK_{\ep1} \setminus \{l_n\},$ where  $\{l_n\}$ is the set
of the first indices of  pairs from $\cN_\ep.$
Then $\{ \la_k:k\in \bK_{\ep 1}\}$ is separated, and so by  Proposition \ref{prop2}, there exists 
 $\{\theta_k, k \in \bK_{\ep1}\},$ a 
family biorthogonal to $\{\exp{i\la_kt} ,\, k \in \bK_{\ep1}\}, $ in the space $L^2(0,T),$ and the following estimate holds
\beq \label{biest3}
\| \theta_k\|_{L^2(0,T)}\leq C_5\exp(C_6(\Im \la_k)^{\varkappa}),
\eeq
with $C_5,C_6$ positive constants that depend only on $T, \al, \rho$.
Now we construct the family $\{h_k(x,t), k \in \bK \}$
in the following way. For $k \in \bK_{\ep} $ we set $$h_k(x,t)=\phi_k(x)\theta_k(t).$$
For $k\notin \bK_{\ep}$, there exist $n,l_n, m_n$ such that 
 $(l_n,m_n) \in \cN_\ep $ and either $l_n=k$ or $m_n=k$. 
 Assume the latter; the argument in the other case is similar. By Lemma \ref{trig}, there exist a pair $\eta_{l_n}(x),\eta_{m_n}(x)$  of functions biorthogonal to $\phi_{l_n}(x),\phi_{m_n}(x)$ on $L^2(a,b)$, and furthermore there exists a positive constant $C$, independent of $n$, such that
\beq
 \int_{a}^b |\eta_{l_n}(x)|^2+|\eta_{m_n}(x)|^2\leq C.\label{biest2}
\eeq
 Let
$$ h_{l_n}(x,t)=\eta_{l_n}(x) \theta_{l_n}(t), \ h_{m_n}(x,t)=\eta_{m_n}(x) \theta_{l_n}(t).$$
It is then easy to check that 
$$\int_0^T \int_a^b h_j(x,t)\exp(i\la_kt)\phi_k(x)\ dxdt=\delta_{jk},\ \forall j,k\in \bK.
$$
Finally, by \rq{biest3} and \rq{biest2},
the estimate \rq{biest1} follows. $\Box$

\ 

We now complete the proof of Theorem \ref{thm2}. Assume for the moment $\la_n^+\neq \la_n^-$ for any $n$.
Recall for $\al <1, \rho<2$ we have $\Im \la_k\asymp |k|^{2\al}$,
while for $\al\in [1,3/2)$, we have 
$$
\Im \la_k= \frac{|k|^{4-2\al}}{\rho}\Big(1+o(1)\Big),\ k>0,\mbox{ and }\Im \la_k= \rho |k|^{2\al}\Big(1+o(1)\Big),\ k<0 .
$$
The moment problem \rq{mo51i} is formally solved by 
$$\sum_{k\in \bK}\zeta_kh_k(x,t).
$$
It suffices to prove convergence of this series. By \rq{mo3i},\rq{mo4i},\rq{tf2i}, there exist positive constants $C_3,C_4$ such that  
$$|\zeta_k|\leq C_3e^{-C_4(\Im \la_k)}, \forall k\in \bK.
$$
Since $\varkappa <1,$
by \rq{biest1} and the asymptotics of $\{ \Im \la_k\}$, the series converges in $L^2((a,b)\times (0,T))$. 

Finally, suppose $\la_n^+=\la_n^-$
for some $n$. It is easy to see this $n$ will be unique.  Here, we need to replace the pair
$(e^{\la_n^+t},e^{\la_n^-t})$ by the pair
$(e^{\la_n^+t},te^{\la_n^+t})$. The construction of the biorthogonal set $\{ \theta_k(t)\}$ can now proceed same as in the previous section, and then the construction of $\{ h_k(x,t)\}$ can now proceed as above in this section. The details are left to the reader. This finishes the proof of Theorem \ref{thm2}.$\Box$

\section{Conclusion}
Our results on finite dimensional control (Theorem \ref{2d}) are in some sense definitive. One possible extension would be to consider perturbations of the Laplacian, replacing $ u_{xx}$ by $(r(x)u_x)_x+q(x)u(x)$.
Another possible extension would be to replace Dirichlet boundary conditions by Neumann or Robin boundary conditions. In all these cases, provided the Sturm-Liouville problem is regular, the spectrum would remain simple,  with the same asymptotics as in this paper. Thus the frequency set will have multiplicity at most two, and two dimensional null-controllabilty will always be possible for $\al <3/2$, and will fail for $\al \geq 3/2.$

In the case of controls distributed on an open proper subset of $(0,\pi)$, Theorem \ref{thm2}, the case $\al \geq 3/2$ remains open. In case the associated Sturm-Liouville problem is perturbed regularly, as in the previous problem, the methods of this  paper will apply. For $\al \in (1,3/2)$,
one would need to generalize Lemma \ref{trig}.

\ 

{\bf Acknowledgements.} The research of Sergei Avdonin was supported in part by the
National Science Foundation, grants DMS 1909869 and 2308377.

\section{Appendix}
{\bf Proof of Lemma \ref{trig}:}
Denote by $\Phi(n,m)$   the angle between  $\f_n$ and $\f_m$ in $L^2(a,b)$, 
\jg{$\Phi(n,m)\in [0,\pi/2]$.}
Evidently, for $m\ne n$, the functions $\sin mx$ and $\sin nx$ are linearly independent on $(a,b)$, which implies $\Phi(n,m)>0$.  Therefore we can restrict to the large $m$, $n$. In what follows we suppose $m>n$.
By the definition 
$$
\cos \Phi(n,m)=\frac{|(\f_n,\f_m)_{L^2(a,b)}|}{\|\f_n\|_{L^2(a,b)}\, \|\f_m\|_{L^2(a,b)}}.
$$
Further,
$$
\int_a^b \sin mx\, \sin nx\, dx=\frac12\left[ 
\frac{\sin (n-m) x}{n-m}- \frac{sin (n+m) x}{n+m}
\right]_a^b.
$$
This gives the asymptotic relation
$$
(\f_n,\f_m)_{L^2(a,b)}=
\frac 2{\pi(m-n)}\sin \left[\frac12 (m-n)(b-a)\right]  \cos \left[ \frac12 (m-n)(b+a)\right]   +O(1/(m+n)).
$$

Similarly
$$
\|\f_m\|_{L^2(a,b)}^2=\frac2\pi\int_a^b \sin^2 mx\, dx=\frac{b-a}\pi+O(1/m).
$$
Then 
$$
\cos \Phi(n,m)=\frac{2}{(b-a)(m-n)}\sin \left[\frac12 (m-n)(b-a)\right]  \cos \left[ \frac12 (m-n)(b+a)\right] 
+O(1/n)\le
$$
$$
\le 
\frac{2}{(b-a)(m-n)}\sin \left[\frac12 (m-n)(b-a)\right] +O(1/n)
$$
The function 
$$
f(x)=\frac{2\sin (x/2)}{x}, 
$$
defined on a semiaxis $(\ep,\infty)$ with a positive $\ep$ satisfies 
$$
\sup f<1.
$$
Indeed, $f(x)<1$  and this functions goes to zero as $x$ goes to infinity,

Thus, 
$$
\sup_{m\ne n} \cos \Phi(n,m)<1.
$$

$\Box$

\end{document}